\documentclass[a4paper,11pt,mathcal]{article}
\usepackage{amsfonts,eucal,amssymb,amsmath,amsthm,amscd} %\input{epsf}

\newtheorem{theorem}{Theorem}[section]
\newtheorem{corollary}[theorem]{Corollary}

\newtheorem{proposition}[theorem]{Proposition}
\newtheorem{lemma}[theorem]{Lemma}
\newtheorem{definition}[theorem]{Definition}

\newenvironment{Pf}{\medskip \noindent {\bf Proof: }}
   {$\diamondsuit$ }

\newcommand{\ad}{\text{ad }}
\newcommand{\Hom}{\text{Hom}}
\newcommand{\End}{\text{End}}
\newcommand{\Mon}{\text{Mon}}
\newcommand{\Ad}{\text{Ad }}
\newcommand{\Aut}{\text{Aut }}
\newcommand{\Autloc}{\ensuremath{\text{Aut}_{loc} }}

\newcommand{\rk}{\text{rk }}
\newcommand{\rkalg}{\text{rk}^{\text{alg}}}
\newcommand{\AdgP}{\text{Ad}_{\mathfrak{g}} P}

\newcommand{\Hol}{\text{Hol}}

\newcommand{\GL}{\text{GL}}

\newcommand{\liep}{\ensuremath{\mathfrak{p}}}

\newcommand{\lieg}{\ensuremath{\mathfrak{g}}}
\newcommand{\lien}{\ensuremath{\mathfrak{n}}}
\newcommand{\liel}{\ensuremath{\mathfrak{l}}}
\newcommand{\liea}{\ensuremath{\mathfrak{a}}}
\newcommand{\liem}{\ensuremath{\mathfrak{m}}}
\newcommand{\liez}{\ensuremath{\mathfrak{z}}}
\newcommand{\lieh}{\ensuremath{\mathfrak{h}}}
\newcommand{\lies}{\ensuremath{\mathfrak{s}}}
\newcommand{\lier}{\ensuremath{\mathfrak{r}}}
\newcommand{\lieu}{\ensuremath{\mathfrak{u}}}

\newcommand{\RR}{\ensuremath{{\bf R}}}
\newcommand{\CC}{\ensuremath{{\bf C}}}
\newcommand{\NN}{\ensuremath{{\bf N}}}
\newcommand{\XX}{\ensuremath{{\bf X}}}
\newcommand{\ZZ}{\ensuremath{{\bf Z}}}

\newcommand{\Had}{\ensuremath{\bar{H}^{\Ad}}}
\newcommand{\Sad}{\ensuremath{\bar{S}^{\Ad}}}

\newcommand{\Sadd}{\ensuremath{\bar{S}_d}}
\newcommand{\Hadd}{\ensuremath{\bar{H}_d}}

\newcommand{\ddtzero}{\ensuremath{\left.\frac{\mbox{d}}{\mbox{d}t} \right|_0}}

\begin{document}
\pagenumbering{arabic}

\title{An embedding theorem for automorphism groups of Cartan geometries}
\author{Uri Bader, Charles Frances, and Karin Melnick}
\date{\today}
\maketitle

\section{Introduction}
The main motivation of this paper arises from classical questions about
actions of Lie groups preserving a geometric structure: which Lie groups can
act on a manifold preserving a given structure, and which cannot?  Which
algebraic properties of the acting group have strong implications on the
geometry of the manifold, such as local homogeneity, or on the dynamics of the
action?  These questions have been studied extensively by many authors (
\cite{AS1}, \cite{AS2},  \cite{DA}, \cite{DAG}, \cite{Gromov},
\cite{kowalsky}, \cite{zeghib1}, \cite{zeghib2}, \cite{Zirigid}, \cite{ZiLor}
among others). These papers primarily prove results in the setting of
either rigid geometric structures of algebraic type, which were defined by M. Gromov
in \cite{Gromov}, or of $Q$-structures of finite type, which were introduced
by Cartan.
% and formalized by Kobayashi \cite{Ko}. 

Our aim in this paper is to develop new tools for studying the kind of
questions mentioned above, in the setting of Cartan geometries. This notion
of geometric structures was developed by E. Cartan in the first half of the last century. A Cartan
geometry infinitesimally models a manifold $M$ on a homogeneous space
$G/P$, where $G$ is a Lie group and $P$ is a closed subgroup.  The geometry is
{\it flat} if and only if $M$ is locally  modeled on $G/P$---that is, $M$ is a
$(G,G/P)$-manifold.  The framework of Cartan geometries offers a very
satisfactory degree of generality, since it essentially comprises all classical rigid geometric structures, including pseudo-Riemannian metrics, conformal pseudo-Riemannian structures, non-degenerate CR structures, and projective structures.  

The main results of this paper are an embedding theorem (theorem
\ref{main-theorem0}) comparable to Zimmer's
well-known embedding theorem for automorphism groups of $Q$-structures; its algebraic consequences on automorphism
groups (theorem \ref{bounds-theorem}); its consequences for
parabolic geometries admitting automorphism groups of high rank (theorem \ref{parabolic-flatness-theorem}); and its
consequences for geometries additionally equipped with a finite continuous volume form (theorems \ref{local-freeness-theorem} and \ref{homogeneous-theorem}). We now discuss in detail these results.

\subsection{Statement and discussion of main theorem}
\label{intro-mainthm-section}

The Lie algebra of a Lie group $G$ is denoted $\lieg$. 

\begin{definition}
\label{def-cartan}
Let $G$ be a Lie group and $P$ a closed subgroup.  A \emph{Cartan geometry} $(M,B,\omega)$ modeled on $(\lieg,P)$ is
\begin{enumerate}
\item{a principal $P$-bundle $\pi: B \rightarrow M$}
\item{a $\lieg$-valued $1$-form $\omega$ on $B$ satisfying}
\begin{itemize}
\item{for all $b \in B$, the restriction $\omega_b : T_bB \rightarrow \lieg$ is an isomorphism} 
\item{for all $b \in B$ and $Y \in \liep$, the evaluation $\omega_b(\left. \frac{\mbox{d}}{\mbox{d} t} \right|_0 be^{tY} ) = Y$ }
\item{for all $b \in B$ and $p \in P$, the pullback $R_p^* \omega = \Ad p^{-1} \circ \omega$}
\end{itemize}
\end{enumerate}
\end{definition}

The automorphism group $\mbox{Aut}(M,B,\omega)$ consists of the bundle
automorphisms of $B$ preserving $\omega$.  The diffeomorphisms of $M$ that
lift to elements of $\mbox{Aut}(M,B,\omega)$ will be denoted $\Aut M$. When the
action of $G$ on $G/P$ has finite kernel, a very mild assumption that we
  will always make in this paper, it can be shown that the projection from
  $\mbox{Aut}(M,B,\omega)$ to $\Aut M$ also has finite kernel. The automorphisms of a framing on a manifold form a Lie
group (see \cite{Ko} I.3.2), from which it follows that $\Aut M$ is a Lie group.

Cartan geometries are closely related, albeit in a nuanced way, to
\emph{$Q$-structures of finite type}.  Kobayashi proved that for $Q$ of \emph{finite type}, the automorphism group of
any $Q$-structure on any $M$ is a Lie group.  For such structures, he
constructs a principal subbundle $B$ of the frame bundle of some order,
equipped with a framing preserved by all automorphisms of the $Q$-structure.
His construction thus associates to a $Q$-structure of finite type a Cartan
geometry with the same automorphism group, although the associated Cartan
geometry is not canonical.  Our theorems on
automorphism groups of Cartan geometries thus apply also to automorphism
groups of $Q$-structures of finite type.

For many differential-geometric structures, including pseudo-Riemannian
metrics, conformal structures, and, more generally, \emph{parabolic
  structures} \cite{CapSch}, there is a canonical associated Cartan geometry
modeled on some $G/P$, and in fact a bijective correspondence between Cartan
geometries modeled on $G/P$ and such structures.  In this case, the
automorphism group of the structure coincides with the automorphism group of
the associated $(M,B,\omega)$. The problem of associating a
canonical Cartan geometry to a geometric structure is called the
\emph{equivalence problem} (see \cite{Sharpe}).  Many interesting examples can
be found in \cite{CapSch}, \cite{chern}, \cite{Sharpe}, \cite{tanaka}.  Some
basic examples are given in section  \ref{intro-parabolic-section}.

\textsc{{\bf Conventions and notations:}}

The Lie algebra of a Lie group will be denoted by the same symbol printed in
the fraktur font; for example, the Lie algebra of $G$ is denoted $\lieg$.

 Throughout the article, $M$ is a connected manifold, and $(M,B,\omega)$ is a
 Cartan geometry modeled on $\XX=G/P$.  We will make the following global assumptions.
\begin{itemize}
\item{\emph{The pair $(G,P)$ is such that the $G$-action on $G/P$ has finite kernel; in practice this assumption is not restrictive.}} 
\item{\emph{The image of $P$ under $\mbox{Ad}$, denoted $\AdgP$, is Zariski closed in $\Aut \lieg$}.}
\item{\emph{For a Lie subgroup $H < \Aut M$, the kernel of $\Ad: H \to
      \Aut \lieh$ is amenable; note that for $H$ connected, $\ker \text{Ad} =
 Z(H)$ is always amenable. }} 
\end{itemize}

 For a subgroup $S < H$, the group $\Sad$ is the Zariski closure of $\Ad S$ in $\Aut \lieh$, and  $\Sadd$ is the maximal Zariski closed
 subgroup of $\Sad$ with no proper cocompact algebraic normal subgroups; see section \ref{background-measure-section} below for details.  The main theorem of the paper is the following.

\begin{theorem}
\label{main-theorem0}
Let $H$ be a Lie subgroup of $\Aut M$.  Let $S$ be a subgroup of $H$ preserving a finite measure on $M$.  Then there are 
\begin{itemize}
\item{$\check{S} < \AdgP$ an algebraic subgroup and} 
\item{$\iota: \lieh \rightarrow \lieg$ a linear monomorphism}
\end{itemize}
such that $\iota(\lieh)$ is $ \check{S}$-invariant, and $\iota$ intertwines $\left.\check{S} \right|_{\iota(\lieh)}$ and $ \Sadd$.
\end{theorem}

Section \ref{main-proof-section} contains a more detailed statement, in theorem \ref{main-theorem}. 

Theorem \ref{main-theorem0} is similar to the embedding theorem of Zimmer
\cite{ZiLor}, which says that when a connected simple group $H$ acts on a
manifold $M$ preserving a finite volume and a $Q$-structure, with $Q$
algebraic, then the standard representation of $Q$ contains the adjoint
representation of $H$.  In particular, $Q$ contains a subgroup locally isomorphic to
$H$. Zimmer's result is proved using algebraic hulls of measurable cocycles
and the Borel density theorem.  An important first step in the proof is to
show that a finite-measure-preserving action of a simple group has discrete
stabilizers almost everywhere.  In this case, the Lie algebra $\lieh$ linearly
embeds in the tangent space at a point, on which $Q$ acts via its standard representation.  

Our result here does not require that $H$ be simple, nor be connected, nor
preserve a finite measure on $M$, nor that the action have any discrete
stabilizers.  The subgroup $S < H$ is not even required to be a Lie
subgroup.  The Cartan connection gives linear embeddings of $\lieh$ in $\lieg$ at every point of $B$.  Then measure-theoretic and algebraic
arguments, including the Borel density theorem, give the existence of
$\check{S} < \AdgP$ restricting to $\Sadd$ on some image of $\lieh$.
Unfortunately, the theorem gives no information for $H$ discrete.

Of course, theorem \ref{main-theorem0} has nice consequences when $H$
preserves a finite measure (see theorems \ref{local-freeness-theorem},
\ref{homogeneous-theorem}, \ref{gromov-corollary} and
\ref{homogeneous-2-theorem}).  Only the subgroup $S$, however, is required to
be measure-preserving.  When $M$ is compact, for example, the theorem applies
to every amenable subgroup $S < H$.  In this way, the theorem yields
information on automorphisms of conformal structures, projective structures,
and other parabolic structures that do not determine a volume form (see theorems \ref{bounds-theorem} and  \ref{parabolic-flatness-theorem}).  

When $H$ is a connected, simple  subgroup of $\Aut M$, and $H$ preserves a
volume on $M$ and a $Q$-structure of finite type, then we can obtain the
conclusion of Zimmer's theorem from theorem \ref{main-theorem0}.
On the other hand, Zimmer's embedding theorem deals with arbitrary
$Q$-structures with $Q$ algebraic; in particular, it does not require finite
type, so there may not be any associated Cartan geometries. With respect to
this comparison, it could be said that we have exchanged the rigidity of simple group actions for rigidity of the structure supporting the action.

\subsection{Bounds on automorphism groups}
\label{intro-bounds-section}

Here we deduce from the main theorem coarse information on the possible groups of automorphisms of $(M,B,\omega)$ when $M$ is compact.  

For a linear group $R < \GL(\RR^n)$, we will define the \emph{real rank}
(resp. the \emph{algebraic real rank}),
denoted $\rk R$ (resp. $\rkalg R$), to be the maximal dimension of a
simultaneously $\RR$-diagonalizable subgroup of $R$ (resp. of the Zariski
closure of a simultaneously $\RR$-diagonalizable subgroup of $R$). The
inequality $\rk R \leq \rkalg R$ always holds.  It can be strict, for example,
for $R \cong \Gamma \ltimes \RR^n$, where $\Gamma$ is a lattice in a semisimple group $G < GL(\RR^n)$ with $\rk G > 0$ (see \cite{PR}, \cite{mostowpr}).
When $R$ is
semisimple, then $\mbox{rk}(\Ad R)=\rkalg (\Ad R)$ is the usual $\RR$-rank of $R$; in fact, for any faithful
representation $\rho$ of $R$, the real rank $\rk \rho(R) = \rkalg \rho(R)$, and
 agrees with the usual $\RR$-rank of $R$.  For $R$ an abstract semisimple group, we will
sometimes write $\rk R$ below without reference to a specific representation. 

 For a Lie group $R$, we define the \emph{nilpotence degree} $n(R)$
to be the maximum nilpotence degree of a connected nilpotent subgroup of $R$.  

\begin{theorem}
\label{bounds-theorem}
Let $(M,B,\omega)$ be a Cartan geometry modeled on $G/P$ with $M$ compact.  Let $H$ be a Lie subgroup of $\Aut M$.  Then
\enumerate
\item{$\rkalg(\Ad H) \leq \mbox{rk}(\mbox{Ad}_\lieg P)$.}  
\item{$n(\Ad H) \leq n(\mbox{Ad}_{\lieg}P)$; in particular, if $H$ is connected, then 
$$n(H)   \leq n(P)+1$$}
\end{theorem}

The bound in (1) is comparable to a theorem of Zimmer for
automorphisms of $Q$-structures.  He proved in \cite{ZiConf} that, for $M$ compact, equipped
with a $Q$-structure, where $Q < \GL(\RR^n)$ is algebraic, any \emph{simple} $H
\subset \Aut M$ has $\rk H \leq \rk Q$.  Note we assume $H$ is neither simple nor connected.  The bound $(1)$ is tight: if $G$ is a simple noncompact Lie group, and $P$ is a parabolic subgroup of $G$, then $\XX=G/P$ is a compact Cartan geometry modeled on itself, on which $G$ acts by automorphisms, and  $(1)$ is an equality in this case. 

The bound in (2) is also tight; an example is found among actions of connected
nilpotent groups by isometries of a compact Lorentz manifold.  Indeed, in this
case, $Q=O(1,n-1)$, the nilpotent subgroups of which are all abelian. The
bound given by the theorem in this case is $n(H) \leq 2$, which is tight since
there exist isometric actions of the Heisenberg group on compact Lorentz
manifolds.  The bound (2) was first proved for Lorentz-isometric actions by Zimmer in \cite{ZiLor}.  For parabolic Cartan geometries, we conjecture the sharper bound $n(H) \leq n(P)$, and we have proved it in the case of conformal pseudo-Riemannian structures (\cite{fm}).

The above theorem does not require that $H$ be connected. Nevertheless, when
$H$ is a discrete group, it yields no information.  The techniques of the
proof do yield results for semi-discrete groups:

\begin{theorem}
\label{bound-semidiscrete-corollary}
Let $L < \GL(V)$ be an algebraic group with no compact algebraic quotients.  Let $\Gamma$ be a
lattice in $L$, and $H = \Gamma \ltimes V$.  Suppose that $H < \Aut M$ for $M$ a compact Cartan geometry modeled on $G/P$.  Then
\begin{enumerate}
\item{$ \rk L \leq \mbox{rk}(\AdgP)$}
\item{$n(L) \leq n(\AdgP)$}
\end{enumerate}
\end{theorem}

\subsection{A rigidity result for parabolic geometries}
\label{intro-parabolic-section}

One class of interesting classical geometries for which the equivalence
problem has been solved is that of {\it parabolic geometries} (see for example
\cite{CapSch}). A Cartan geometry $(M,B,\omega)$ modeled on $\XX=G/P$ is said to be {\it parabolic} when $G$ is a real semisimple Lie group and $P$ is a parabolic subgroup of $G$.  The main examples of such geometries follow.

$\bullet${ \it Conformal pseudo-Riemannian structures}. 

\indent A conformal pseudo-Riemannian structure of signature $(p,q)$ on a
manifold $M$ is the data of a conformal class of metrics 
$$[g]=\{ e^{\sigma}g \ | \ \sigma \in C^{\infty}(M)  \}$$

where $g$ is a metric of signature $(p,q)$. These structures are in canonical
correspondence with Cartan geometries modeled on $O(p+1,q+1)/P$, where 
$$P \cong (\RR* \times O(p,q))\ltimes {\RR}^{p+q}$$

is a maximal parabolic subgroup of $O(p+1,q+1)$, and the connection $\omega$
satisfies certain normalization conditions.

$\bullet$ {\it Nondegenerate $CR$-structures}.

These structures model real hypersurfaces in complex manifolds.  A
nondegenerate $CR$-structure of type $(p,q)$ on a $(2m+1)$-dimensional
manifold $M$, where $p+q=m$, is the data of a contact subbundle $E \subset
TM$ equipped with an almost-complex structure $J$ and a conformal class of
Hermitian metrics of type $(p,q)$.  Such a structure is equivalent to a unique
normal Cartan geometry modeled on $SU(p+1,q+1)/P$, where $P$ is the parabolic
subgroup of $SU(p+1,q+1)$ stabilizing an isotropic complex line in ${\bf
  C}^{p+1,q+1}$.  Again, the connection $\omega$ must satisfy some
normalizaiton conditions.  The equivalence problem for strictly pseudo-convex
$CR$-structures was first solved by E. Cartan in dimension $3$, and in the general case in \cite{chern}, \cite{tanaka}, \cite{CapSch}. 

$\bullet$ {\it Projective structures.} 

A \emph{projective structure} on an $n$-dimensional manifold $M$ is a family of smooth curves $C:I \to M$, defined locally by an ODE of the form
\[ \frac{c_1^{\prime \prime}+P_1(c^{\prime})}{c_1^{\prime}}= \frac{c_2^{\prime
    \prime}+P_2(c^{\prime})}{c_2^{\prime}}=...= \frac{c^{\prime \prime}_n +P_n(c^{\prime})}{c_n^{\prime}} \]

with $P_k(x)=\Sigma a_{ij}^k x_ix_j$ for smooth functions $a_{ij}^k$.
Geodesics of a Riemannian or pseudo-Riemannian metric, for example, satisfy
such an equation. To such a projective structure corresponds a unique normalized
Cartan geometry modeled on the projective space ${\RR}P^n=PGL(n+1,{\RR})/P$,
where $P$ is the stabilizer of a line in $\RR^{n+1}$ (see \cite{Sharpe} ch. $8$ for details).   
 
We now come to our main result concerning parabolic geometries. Note that such
a geometry $(M,B,\omega)$ is modeled on $G/P$ with $\mbox{rk}(\AdgP) = \rk G$.
Theorem
\ref{bounds-theorem} asserts that when a Lie group $H$ acts by
automorphisms of $(M,B,\omega)$, then
$\rkalg(\Ad H) \leq \rk G$.  When $H$ is connected, and equality
$\mbox{rk}(\Ad H) = \rk G$ holds, a strong rigidity phenomenon occurs.  The universal covers of $\XX$ and $G$ are denoted
$\widetilde{\XX}$ and $\widetilde{G}$, respectively.

\begin{theorem}
\label{parabolic-flatness-theorem}
Let $(M,B,\omega)$ be a compact parabolic Cartan geometry, modeled on
$\XX=G/P$, where $G$ has finite center, and $\omega$ is regular.  Let $H< \Aut M$ be a Lie subgroup. 
\enumerate
\item{If $\rkalg (\Ad H) = \rk G $, then there is a Lie algebra embedding $\lieh \to \lieg$.}
\item{If $H$ is connected, and  $\mbox{rk}(\Ad H) = \rk G$, then $M \cong
    \Gamma\backslash \widetilde{\XX}$, for some discrete subgroup $\Gamma < \widetilde G$.} 
\end{theorem}

The regularity condition on the Cartan connection is very mild. For all
parabolic geometries for which the equivalence problem has been solved,
regularity is one of the normalization conditions
on the Cartan connection.  See section
\ref{parabolic-background-subsection} for the definition of regular
connection.

Theorem \ref{parabolic-flatness-theorem} is a wide generalization of previous
results obtained in \cite{BN} and \cite{fz} for conformal actions of connected simple groups on compact
pseudo-Riemannian manifolds.  For parabolic Cartan geometries modeled on $G/P$ with $\rk
G = 1$, existence of a nonproper automorphism group implies that $M$ is equivalent
to $G/P$ or the complement of a point in $G/P$.  In conformal
Riemannian geometry, this statement is a celebrated theorem of Lelong-Ferrand
\cite{LF}; it is proved for general rank-one parabolic Cartan geometries in \cite{frances}.

\subsection{Automorphism groups preserving a finite volume}
\label{intro-unimodular-section}

By a finite volume on a manifold $M$, we mean a continuous volume form $\nu$
on $M$, such that $\int_{M}\nu < \infty$. This section contains results on
actions by automorphisms of a Cartan geometry preserving a
finite volume (Theorem \ref{local-freeness-theorem} makes the slightly weaker assumption that the action preserves a finite measure of full support).  Isometric actions on pseudo-Riemannian manifolds
preserve a volume form, which is always finite if the manifold is compact.
Other interesting examples of such actions are affine actions preserving a
finite volume.  Note that the invariant volume form need only be $C^0$. 

Recall that $\Hadd$ is the discompact radical of $\Had$ (see \S
\ref{background-measure-section}).  In \cite{DAG}, it is conjectured that any smooth, faithful, volume-preserving action of a simple Lie group on a compact manifold is everywhere locally free. Here, we prove:

\begin{theorem}
\label{local-freeness-theorem}
Let $(M,B,\omega)$ be a  Cartan geometry, and let $H < \Aut M$ be a Lie subgroup preserving a finite measure of full support in $M$.  Then for almost-every $x \in M$, the
stabilizer Lie algebra $\lieh_x$ is $\Hadd$-invariant.  If
there are only finitely-many $\Hadd$-invariant subalgebras of $\lieh$, then the
$H$-action is locally free on a dense, open, conull subset of $M$.
\end{theorem}

This theorem is proved in \cite{Zifree} for connected semisimple Lie groups
with no compact factor.  Our conclusion of local freeness here holds for a wider class of groups, including semidirect products $H=\Gamma \ltimes_\psi R$, where $R$
is a connected Lie group, $\Gamma$ is a lattice in a semisimple group $L$, and
$\psi$ is the restriction of a nontrivial representation of $L$ in $\Aut \lier$.

The next theorem says that when a group of automorphisms of a Cartan geometry preserves a finite volume, then existence of an open orbit implies homogeneity; in this case, the action is
everywhere locally free.  In general, it is an interesting open question for a
geometric manifold when existence of an open locally homogeneous subspace
implies local homogeneity everywhere (\cite{DAG}, \cite{Du}, \cite{Me}).

Recall that $\Had$ denotes the Zariski closure of $\Ad H$, and $\Hadd$ its
discompact radical. 

\begin{theorem}
\label{homogeneous-theorem}
Let $(M,B,\omega)$ be a  Cartan geometry, and let
$H < \Aut M$ be a Lie subgroup such that $\Had = \Hadd$.  If $H$ preserves a
finite volume on $M$ and has an open orbit on $M$, then
\enumerate
\item{ The action of $H^o$, the identity component of $H$, is transitive on $M$. }
%Moreover $H^o=(\Aut M)^o$.}
\item{  The action of $H$ is everywhere locally free.}
\item{The stabilizers of points in $H^o$ are lattices  $\Gamma \subset H^o$.} 
\end{theorem}

Essentially all classical rigid geometric structures can be described both as
Cartan geometries and rigid geometric structures of algebraic type in the
sense of Gromov.  The interplay of these two points of
view can be fruitful, as the following corollary shows.  

\begin{corollary}
\label{gromov-corollary}
Let $(M,B,\omega)$ be a $C^\omega$ Cartan geometry with $M$ simply
connected. Assume that there is a natural associated rigid structure of algebraic type $\cal{S}$
on $M$ with $\Autloc M =\Autloc (M,\cal{S})$. Let $\nu$ be a $C^\omega$ finite
volume form on $M$.  Let $H < \Aut M$ be a Lie subgroup such that $\Had = \Hadd$.  If $H$ preserves $\nu$ and has a dense orbit
in $M$, then the geometry is homogeneous: there exists $H^{\prime} < \Aut M$ acting transitively on $M$ and preserving $\nu$. 
\end{corollary}

\subsection{Homogeneous reductive geometries}

Finally, we record a corollary for those Cartan geometries called reductive. A Cartan geometry $(M,B,\omega)$ modeled on $G/P$ is \emph{reductive}
when there is an $\AdgP$-invariant decomposition $\lieg=\lien \oplus \liep$, in which the restriction of $\AdgP$ to $\lien$ is faithful.  Note that $\lien$ is in general a subspace of $\lieg$, not a subalgebra.
% When considering reductive geometries, we will always assume that the action
% of $P$ on $\lien$ by the adjoint representation is faithful. 
Classical examples of reductive geometries are pseudo-Riemannian metrics, for
which 
$$(\lieg,\liep)=(\mathfrak{so} (p,q) \ltimes \RR^n,\mathfrak{so}(p,q))$$
or linear connections, for which 
$$(\lieg,\liep)=(\mathfrak{gl}(\RR^n) \ltimes \RR^n,\mathfrak{gl}(n,\RR))$$

\begin{corollary}
\label{homogeneous-2-theorem}
Let $(M,B,\omega)$ be a  reductive Cartan geometry, and $\lieg=\lien \oplus
\liep$ an $\AdgP$-invariant decomposition.  Let $H < \Aut M$ be a Lie subgroup
such that $\Had = \Hadd$.  If $H$ acts transitively on $M$ preserving a finite
volume, then
\enumerate
\item{The group $\Had$ is isomorphic to a subgroup $ \check H <  P$.}
\item{There is a linear isomorphism $\iota : \lieh \to \lien$ with $\Ad \check
    H$-invariant image, which intertwines the representations of $\Had$ on $\lieh$ and of $\check H$ on $\iota(\lieh)$}.
\item{ The identity component $H^o$ is isomorphic to a central extension of a subgroup $H^{\prime} < \check H$}.
\end{corollary}

\section{Automorphisms of a Cartan geometry}
\label{background-cartan-section}
Let $(M,B,\omega)$ be a Cartan geometry modeled on $\XX=G/P$, and $H < \Aut M$
a Lie subgroup.  The $H$-action on $B$ is free because it preserves a framing, given by the pullback by $\omega$ of any basis of $\lieg$ (see \cite{Ko} I.3.2).

This section introduces several important $H$-equivariant maps from $B$ to
algebraic varieties.  These maps will be central in the proof of the main theorem.

\subsection{The curvature map of a Cartan geometry}
 The simplest example of a Cartan geometry is the triple $(\XX,G,\omega^G)$,
 where $\omega^G$ is the Maurer-Cartan form on $G$, which evaluates to $X$ on 
the left-invariant vector field corresponding to $X \in \lieg$.  In addition
to fulfilling the definition of a Cartan geometry, the Maurer-Cartan form on $G$ satisfies the {\it structural equation}:

\[ d\omega^G(X,Y)+[\omega^G(X),\omega^G(Y)]=0  \] 

for every pair of vector fields $X$ and $Y$ on $G$.

For a general Cartan geometry $(M,B,\omega)$, the 2-form 
$$K(X,Y)=d\omega(X,Y)+[\omega(X),\omega(Y) ]$$

is not zero, and is called {\it the curvature form} of $\omega$.  It is a
fundamental fact that $K$ vanishes if and only if $(M,B,\omega)$ is locally
isomorphic to $(\XX,G,\omega^G)$ (see \cite{Sharpe}). Thus, the manifolds $M$
{\it locally} modeled on $\XX$, also known as manifolds with $(G,\XX)$-structures, correspond exactly to {\it flat} Cartan geometries $(M,B,\omega)$, modeled on $\XX$.  

The associated {\it curvature map} is 
\begin{eqnarray*}
\kappa & : & B \to \Hom(\Lambda^2(\lieg/\liep),\lieg) \\
\kappa_b(\sigma(X), \sigma(Y)) & = & K_b(\omega_b^{-1}(X),\omega_b^{-1}(Y))
\end{eqnarray*}
Where $X,Y \in \lieg$, and $\sigma$ is the projection $\lieg \rightarrow \lieg
/ \liep$.  This
map is well-defined because $K_b$ vanishes as soon as one of its argument is
tangent to a fiber of $B \to M$. The map $\kappa$ is $P$-equivariant, where
$P$ acts on $\Hom(\Lambda^2(\lieg/\liep),\lieg)$ via the adjoint: for any  $\check p \in \Ad P$ and $\beta \in \Hom(\Lambda^2(\lieg/\liep),\lieg)$,

$$(\check p.\beta)(u,v)={\check p}^{-1}.(\beta(\check p.u,\check  p.v))$$

 For any $b \in B$, the Cartan connection induces a natural linear map
\begin{eqnarray*} 
\iota_b & : & \lieh \to \lieg \\ 
\iota_b(X) & = & \omega_b (\ddtzero \varphi_X^t.b)
\end{eqnarray*}

where $\varphi_X^t$ is the flow along $X$ for time $t$.  The freeness of the
$H$-action implies that $\iota_b$ is injective. In general, $\iota_b$ is not a
Lie algebra homomorphism.  It nearly is, however, when $\kappa_b = 0$.  

\begin{lemma}
\label{curvature-lemma}
Let $(M,B,\omega)$ be a Cartan geometry modeled on $G/P$. If  for some $b \in B$, the curvature
$\kappa_b$ vanishes on the image of $\sigma \circ \iota_b$, then $- \iota_b:
\lieh \to \lieg$ is an embedding of Lie algebras.
\end{lemma}

% Recall that an anti-homomorphism $\iota$ of Lie algebras is a linear map of
% the underlying vector spaces satisfying 
% $$ \iota([X,Y]) = [\iota(X),\iota(Y)]$$

% Given an anti-homomorphism $\iota : \lieh \rightarrow \lieg$, the map
% $$ \iota^- : X \rightarrow - \iota(X)$$

% is a homomorphism from $\lieh$ to $\lieg$.

 \begin{Pf}
Let $X,Y \in \lieh$, viewed as Killing fields on $B$.  Recall that 
$$ d\omega(X,Y) = X.\omega(Y) - Y.\omega(X) - \omega([X,Y])$$

Since $X$ is a Killing field for $\omega$,
$$ \omega(Y(b)) = \omega ((\varphi_X^t)_{*b}Y(b))$$

Now 
\begin{eqnarray*}
(X.\omega(Y))(b) & = & \ddtzero \omega(Y(\varphi_X^t b)) \\
%% & = & \ddtzero \left[ \omega((\varphi_X^t)_{*b} Y(b)) + \omega(Y(\varphi_X^t(b))) - \omega((\varphi_X^t)_{*x}Y(b)) \right] \\
& = & \ddtzero \left[ \omega(Y(\varphi_X^t b)) - \omega((\varphi_X^t)_{*b}Y(b)) \right] \\
& = & \omega([X,Y](b))
\end{eqnarray*}

Similarly, 
$$ (Y.\omega(X))(b) = \omega([Y,X](b))$$

Then
\begin{eqnarray*}
\kappa_b(\sigma \circ \iota_b(X),\sigma \circ \iota_b(Y)) & = & d \omega(X(b),Y(b)) + [\omega(X(b)),\omega(Y(b))] \\
& = & \omega([X,Y](b)) - \omega([Y,X](b)) - \omega([X,Y](b)) + [\omega(X(b)),\omega(Y(b))] \\
& = & \omega([X,Y](b)) + [\omega(X(b)),\omega(Y(b))] \\
 & = & \iota_b([X,Y])+ [\iota_b(X),\iota_b(Y)] 
\end{eqnarray*}

Therefore, vanishing of $\kappa_b$ implies
$$ - \iota_b([X,Y])= [ \iota_b(X),\iota_b(Y)]= [- \iota_b(X), - \iota_b(Y)]$$
as desired.
\end{Pf}

% \begin{Pf}
% Let $X,Y \in \lieh$.  Recall that 
% $$ d\omega(X,Y) = X.\omega(Y) - Y.\omega(X) - \omega([X,Y])$$

% Because $X$ is a Killing field for $\omega$,
% $$ \omega(Y(b)) = \omega ((\varphi_X^t)_{*b}Y(b))$$

% Now 
% \begin{eqnarray*}
% (X.\omega(Y))(b) & = & \ddtzero \omega(Y(\varphi_X^t b)) \\
% %% & = & \ddtzero \left[ \omega((\varphi_X^t)_{*b} Y(b)) + \omega(Y(\varphi_X^t(b))) - \omega((\varphi_X^t)_{*x}Y(b)) \right] \\
% & = & \ddtzero \left[ \omega(Y(\varphi_X^t b)) - \omega((\varphi_X^t)_{*b}Y(b)) \right] \\
% & = & \omega([X,Y](b))
% \end{eqnarray*}

% Similarly, 
% $$ (Y.\omega(X))(b) = \omega([Y,X](b))$$

% Then
% \begin{eqnarray*}
% \kappa_b(\sigma \circ \iota_b(X),\sigma \circ \iota_b(Y)) & = & d \omega(X(b),Y(b)) + [\omega(X(b)),\omega(Y(b))] \\
% & = & \omega([X,Y](b)) - \omega([Y,X](b)) - \omega([X,Y](b)) + [\omega(X(b)),\omega(Y(b))] \\
% & = & \omega([X,Y](b)) + [\omega(X(b)),\omega(Y(b))] \\
% & = & \iota_b([X,Y])+[\iota_b(X),\iota_b(Y)]
% \end{eqnarray*}

% Therefore, vanishing of $\kappa_b$ implies
% $$- \iota_b([X,Y])=[- \iota_b(X),- \iota_b(Y)]$$

% as desired.
% \end{Pf}

\subsection{Natural equivariant maps to varieties}
\label{natural-objects-section}
 
To an action of $H$ by automorphisms of a Cartan geometry $M$ modeled on $G/P$
can be associated
various equivariant maps to algebraic varieties built from the adjoint
representations of $H$ and $G$.

The maps $\iota_b:\lieh\to\lieg$ for $b \in B$ give a map
$$\iota: \ B \to \Mon(\lieh,\lieg)$$
 where $\Mon(\lieh,\lieg)$ is the variety of all injective linear maps from $\lieh$ to $\lieg$. The map $\iota$ is continuous, since $\omega$ is so.
Observe that $H$ and $P$ have commuting actions on $B$,
as well as on $\Mon(\lieh,\lieg)$, by pre- and post-compositions with the
adjoint actions.  It will be crucial in the sequel that the latter $(\Ad H
\times \AdgP)$-action is the
 restriction of an algebraic action, and that $\AdgP$ is algebraic.
To a Lie subgroup $H < \Aut M$, can be associated a real-algebraic space 
$$U=\Mon(\lieh,\lieg)\times \Hom(\wedge^2(\lieg/\liep),\lieg)$$
where $H$ acts trivially on $\Hom(\wedge^2(\lieg/\liep),\lieg)$, and a
continuous map $\phi : \ B \to U$, which is continuous and $(H \times
P)$-equivariant.

We denote by $\check{P}^b$ the algebraic subgroup of $\AdgP$ consisting of
those $\check{p}$ such that
$\check{p}.\iota_b(\lieh)=\iota_b(\lieh)$. There is a natural algebraic
homomorphism 
\begin{eqnarray*}
\rho_b & : & \check{P}^b \to \GL(\lieh) \\
\rho_b(\check{p}) & = & \iota_b^{-1} \circ \check{p} \circ \iota_b
\end{eqnarray*}

Now, if $S < H$ is a subgroup, and $b \in B$, define

\[ \check{S}^b = \{ \check p \in \check{P}^b \ | \ \rho_b(\check p) \in \Sadd \ \mbox{and} \  \check p.\kappa_b=\kappa_b \} \]

The group $\check{S}^b$ may not be algebraic in general, but it is if $\rho_b:
\check{S}^b \to \Sadd$ is {\it surjective}, as in this case it is the preimage
in $\check{P}^b$ of the algebraic group $\Sadd$ by the algebraic homomorphism
$\rho_b$.  The restriction of $\rho_b$ to  $\check S_b$ is then an algebraic
epimorphism onto $\Sadd$. 

The main point in the proof of theorem \ref{main-theorem0} will be precisely to show the existence of $b \in B$ such that $\rho_b: \check S_b \to \Sadd$ is  a surjection.

\section{Background on measure theory}
\label{background-measure-section}

In this section we present the version of the Borel density theorem that
figures in the proof of the main theorem \ref{main-theorem0}.
Primary references for material here on invariant measures for real-algebraic actions are \cite[Chapters 2,3]{ZiETSG} and \cite{Shalom}.
We will be considering real-algebraic groups and varieties, and all
results will be stated in this setting; in fact, we will use the term
``algebraic'' to mean real-algebraic below.

For $G$ an algebraic group and $V$ an algebraic variety, an algebraic action
of $G$ on $V$ is an action given by an algebraic morphism $\varphi : G \times V \rightarrow V$. 
For a locally compact group $S$ and a homomorphism $\psi:S\to G$, where $G$ is an algebraic group, 
we denote the real points of the Zariski closure of $\psi(S)$ by $\bar{S}^{\psi}$ or $\bar{S}$.

% By a {\em real algebraic action} we mean a triple $({\bf V},S,\rho)$,
% where ${\bf V}$ is an algebraic variety defined over $\RR$, $S$ is a locally compact group and $\rho$ is a homomorphism
% from $S$ to $(\mbox{\bf Aut V})(\RR)$ - the real point of the algebraic group of automorphisms of ${\bf V}$.
% Clearly, via $\rho$, $S$ acts on $V={\bf V}(\RR)$.
% Abusing the notation we will say that $S$ acts algebraically on $V$. 

The \emph{discompact radical} $G_d$ of an algebraic group $G$ is the maximal
algebraic subgroup of $G$ with no cocompact, algebraic, normal subgroups.  It
is the minimal algebraic subgroup $Q$ containing all algebraic subgroups with no
compact algebraic quotients, which exists by the Noetherian property of
algebraic subgroups.  By convention, we will put $G_d=\{ e \}$ when $G=\{ e \}$.
%% Then any algebraic homomorphism onto a compact group contains $Q$
%% in the kernel, so $Q=G_d$.  
Note that $G_d$ is characteristic.  See
\cite{Shalom} for more details.  Given $\psi : S \rightarrow G$, we
denote by $\Sadd^\psi$ or $\Sadd$ the discompact radical of $\bar{S}^{\psi}$.

The formulations of theorems \ref{main-theorem0} and
\ref{local-freeness-theorem} in terms of $\Hadd$ and $\Sadd$ apply to a wider
class of groups than those for which the Zariski closures $\Had$ or
$\Sad$ have no compact algebraic quotients.  Consider the following
group.  Let 
$$ g = \left( \begin{array}{ccc}
       \cos \theta & - \sin \theta & 0 \\
       \sin \theta & \cos \theta & 0 \\
       0          & 0    &   e^{\lambda}
 \end{array} \right)  \in GL(3,\RR)
$$

for some $\theta \neq 0$ modulo $2\pi$, $\lambda \neq 0$.  Let $H = \langle g \rangle \ltimes \RR^3
\cong \ZZ \ltimes \RR^3$.  Then $\Hadd^{\Ad} \cong \RR^*$, but for any $S < H$, the
group $\Sad$ is either trivial or has $S^1$ as a direct factor.

\begin{theorem}[{\cite[3.11]{Shalom}}] \label{thm:shalom}
Let $\psi : S \rightarrow \Aut V$, for $S$ a locally compact group and $V$ an
algebraic variety.  Assume $S$ preserves a finite measure
$\mu$ on $V$.  Then $\mu$ is supported on the set of $\Sadd$-fixed points in $V$.
\end{theorem}

\begin{Pf}[sketch]
By an ergodic decomposition argument, it is enough to assume that the measure $\mu$ is $S$-ergodic.
It is well known that $\bar{S}$-orbits in $V$ are locally closed \cite[3.1.1]{ZiETSG}, and it follows that
$\mu$ is supported on a single $\bar{S}$-orbit.
We claim that $\mu$ is $\bar{S}$-invariant.
Indeed, Zimmer proved that the stabilizers of probability measures are Zariski
closed when $V$ is a quasi-projective variety \cite[3.2.4]{ZiETSG},
but by a theorem of Chevalley every $\bar{S}$-orbit is $\bar{S}$-equivariantly quasi-projective.
It remains to show that every $\bar{S}$-invariant measure is supported on the set of $\Sadd$-fixed points.
The latter group is generated by algebraic $1$-parameter groups, isomorphic to either the additive or the multiplicative group of $\RR$.
The statement is true for these groups, again by applying an ergodic decomposition and restricting to one orbit.
The proof follows.
\end{Pf}

\begin{corollary} \label{cor:exist-fixed-point}
Let $\psi : S \rightarrow \Aut V$, for $S$ a locally compact group and $V$ an
algebraic variety.
Suppose $S$ acts continuously on a topological space $M$ preserving a finite
Borel measure $\mu$.
Assume $\phi:M\to V$ is an $S$-equivariant measurable map.
Then $\phi(x)$ is fixed by $\Sadd$ for $\mu$-almost-every $x \in M$.
\end{corollary}

\begin{Pf}
Follows from the previous theorem, with the measure $\phi_*(\mu)$.
\end{Pf}

\section{Proof of the main theorem \ref{main-theorem0}}
\label{main-proof-section}

This section will be devoted to the proof of theorem \ref{main-theorem0}, or more precisely to that of theorem~\ref{main-theorem} below which is  a detailed version of it.

We consider $(M,B,\omega)$,  a fixed Cartan geometry 
modeled on $G/P$. Recall the global assumption that $\AdgP
< \Aut \lieg$ is Zariski closed. The notations are those of section
\ref{natural-objects-section}.  Recall that $\pi: B \rightarrow M$ is the bundle projection.

%% KARIN CHANGED WORDING A BIT (NOT CONTENT)
%%%%%%%%%%%%%%%%%%%%%%%%%%%%%%%%%%%%%%%%%%%%
\begin{theorem}
\label{main-theorem}
Let $H < \Aut M$ be a Lie subgroup and $S < H$.
Suppose there is an $S$-invariant measure $\mu$ on $M$.
Then there exists $\Lambda \subset B$ such that $\pi(\Lambda)$ is $\mu$-conull
and, for every $b\in \Lambda$, 
$$\rho_b : \check{S}^b \rightarrow \Sadd$$
is an algebraic epimorphism.  In particular, $\check{S}^b$ preserves $\kappa_b$ and $\rho_b$ intertwines the representation of $\check{S}^b$ on
    $\iota_b(\lieh)$ and that of $\Sadd$ on $\lieh$.
\end{theorem}

%REMOVED BY CHARLES. I PUT IT IN THE THEOREM\begin{remark}
%The proof below actually shows (uppon applying Theorem~\ref{thm:shalom} instead of Corollary~\ref{cor:exist-fixed-point})
%that the upshot of the theorem is valid for every $b$ in the $\pi$-pre-image of $\mu$-almost every $m\in M$.
%The information given by the theorem  is of particular usefulness in the following two cases:

%$\bullet$ If $\mu$ is fully supported---for example, $\mu$ is a volume measure.

%$\bullet$
%If $S$ is amenable.
%Then such $b$ exist in the preimage of every invariant compact subset of $M$, as
%every invariant compact subset of $M$ supports an $S$-invariant measure.
%\end{remark}

\begin{Pf}
Observe first of all that the main point is to prove  the surjectivity of
$\rho_b$ at some point.  The two other assertions of the theorem follow by the definitions of $\check{S}^b$ and $\rho_b$.

As in section \ref{natural-objects-section}, consider the real-algebraic variety 
$$U=\Mon(\lieh,\lieg)\times \Hom(\wedge^2(\lieg/\liep),\lieg)$$
and the $H\times P$-equivariant map
$\phi=\iota \times \kappa:B\to U$.

The map $\phi$ descends to an $H$-equivariant map between the spaces of $P$-orbits:
\[ \bar{\phi}:M=B/P \to V=U/P \]

%% STRATIFICATION MODIFIED BY KARIN
%%%%%%%%%%%%%%%%%%%%%%%%%%%%%%%%%%%
% For $i=0,\ldots,\dim(P)$ denote by $U_i$ the subset of $U$ consisting of points with $i$-dimensional stabilizer in $P$,
% that is
% \[ U \supset U_i=\{ u\in U~|~\dim(\Stab_P(u))=i \}\]
% By the fact that $P<\GL(\lieg)$ is Zariski closed,
% the decomposition $U=U_0 \cup U_1 \cup \cdots \cup U_n$ is a $P$-invariant locally closed stratification.
% The $H$ action commutes with the $P$ action, hence each stratum is $H$-invariant as well.

The $\AdgP$-action on $U$ is algebraic, and there is a $P$-invariant
stratification
$$ U = U_0 \cup \cdots \cup U_r$$

such that, for each $j \geq 0$, the stratum $U_j$ is Zariski open and dense in $\cup_{i\geq j} U_i$, and each quotient space $V_i = U_i / \AdgP$ is a smooth algebraic variety
(see \cite{Ros} or \cite{Gromov} 2.2).  Because the $H$-action on $U$
commutes with $\AdgP$, it preserves the stratification and acts algebraically on each $V_i$.

% Because $\AdgP$ acts algebraically on $U$, the partition by orbit type
% $$ U_Q = \{ x \in U \ : \ g P(x) g^{-1} = Q \ \  \mbox{for some} \ g \in P \}$$ 
% as $Q$ ranges over algebraic subgroups of $\AdgP$, defines a locally finite, hence
% countable, stratification of $U$ (see REF).  For each $Q = Q_i$, the quotient
% $V_i= U_Q/P$ is a smooth algebraic variety.  
% It is standard that, for each $i$, the orbits space $V_i=U_i/P$ has a (real) algebraic structure,
% and an induced  (real) algebraic action of $H$.

Denote by $\mu_i$ the restriction of $\mu$ to $M_i=\bar{\phi}^{-1}(V_i)$.  It
suffices to prove that $\rho_b$ is surjective for $\mu_i$-almost-every $\pi(b)
\in M_i$, for all $i$.
Because $M_i$ is $S$-invariant for all $i$, each measure $\mu_i$ is $S$-invariant.
Given $i$ such that $\mu_i\neq 0$, corollary~\ref{cor:exist-fixed-point} applied to the $S$-equivariant map
$\bar{\phi}:(M_i,\mu_i)\to V_i$ gives that $v=\bar{\phi}(m)$ is an $\Sadd$-fixed point for $\mu_i$-almost-every $m \in M_i$.
For such an $m$, choose $b\in \pi^{-1}(m)\subset B$, and set $u=\phi(b)$.
Then $\AdgP.u$ is stabilized by $\Sadd$. This means that for any $s \in \Sadd$,
there exists $\check{p}_s \in \AdgP$, not necessarily unique, such that $s.(\iota_b,\kappa_b)=\check{p}_s.(\iota_b,\kappa_b)$, or in other words, $\iota_b \circ s=\check p_s \circ \iota_b$, and $\check p_s.\kappa_b=\kappa_b$. Thus, $\check p_s \in \check{S}^b$, and $\rho_b(\check p_s)=s$. We thus get that $\rho_b : \check{S}^b \to \Sadd$ is a surjection, and theorem \ref{main-theorem} is proved.  \end{Pf}

\section{Proofs of bounds on automorphism groups}

Recall that $\Hadd$ denotes the discompact radical of the Zariski closure of
$\Ad H$.  A Lie subalgebra of $\End(V)$ will be called {\it algebraic} when it
is the Lie algebra of a Zariski closed subgroup of $\GL(V)$. Also, throughout this section, $M$ is compact. 

%% LEMMA WRITTEN BY KARIN, REMOVED BY KARIN
%% BECAUSE THM DOESN'T HOLD FOR \Hadd
% The following algebraic lemma establishes the necessary relation between the
% rank and nilpotence degree of $H$ and those of $\Hadd$.

% \begin{lemma}
% \label{bounds-alg-lemma}
% Let $H$ be a Lie group.
% \begin{enumerate}
% \item{$\mbox{rk}(\Ad H) \leq \rk \Hadd$}
% \item{$n(H) \leq n(\Hadd) + 1$}
% \end{enumerate}
% \end{lemma}

% \begin{Pf}
% For (1), let $A$ be a maximal $\RR$-split torus of $\Ad H$, and $\bar{A}$ its Zariski closure.  Now $\bar{A}$ is still abelian and $\RR$-split, so it clearly has not compact algebraic quotients.  Therefore, 
% $$A \subseteq \bar{A} \subseteq \Hadd$$
% yielding the desired inequality.

% For (2), note that if $N$ is nilpotent, then $\Ad N$ is unipotent.  The Zariski closure $\bar{N}$ of $\Ad N$ is still unipotent and clearly has no compact algebraic quotients. Therefore,
% $$ \Ad N \subseteq \bar{N} \subseteq \Hadd$$
% When $N$ is connected, then $(\ker \Ad)\cap N$ is
% contained in $Z(N)$, so $n(N) = n(\Ad N) + 1 \leq n(\Hadd) + 1$, as desired.
% \end{Pf}

\subsection{Proof of theorem \ref{bounds-theorem}}
\label{bounds-proof-section}
To prove (1) of theorem \ref{bounds-theorem}, choose  an $\RR$-split subgroup
of $\Ad H$ with Zariski closure $\bar{S}$ of dimension $ \rkalg(\Ad H)$. Let $S = \mbox{Ad}^{-1}(\bar{S}) < H$.
The kernel of $\mbox{Ad}$ is amenable by assumption (\S \ref{intro-mainthm-section}), and $\bar{S}$ is abelian, so $S$ is amenable.
There is a finite $S$-invariant measure on $M$, and theorem \ref{main-theorem}
gives an algebraic subgroup $\check{S}^b < \AdgP$ such that $\rho_b(\check{S}^b) =
\Sadd = \bar{S}$.  

The Lie algebra of $\check{S}^b$ admits a Levi decomposition into algebraic
subalgebras $\check{\lies}^b=\lier \ltimes \lieu$, where $\lieu$ consists of nilpotent endomorphisms, and $\lier$ is a reductive
 Lie algebra (see \cite{witteratner} ch.$4$).  The homomorphism $\rho_b$ induces an epimorphism of
 algebraic algebras $d\rho_b: \check{\lies}^b \to \bar{\lies}$. The image of
 $\lieu$ consists of elements that are both $\RR$-split and nilpotent, so that
 $\lieu \subset \ker d\rho_b$.  The algebra $\lier$ splits into algebraic
 subalgebras $\lier=\liez \oplus \liel$ where $\liel$ is semisimple, and $\liez$ is an abelian algebra
 of $\CC$-split endomorphisms.  Because $\bar{\lies}$ is abelian, $\liel
 \subset \ker d\rho_b$.  The $\CC$-split algebraic subalgebra $\liez$
 decomposes as a sum $ \liez = \liez^s \oplus \liez^e$ where $\liez^s$ consists of $\RR$-split elements, and $\liez^e$ of elements
 with purely imaginary eigenvalues.  Since $d \rho_b$ is a morphism of real
 algebraic Lie algebras, it respects this decomposition.  Therefore, $\liez^e
 \subset \ker d\rho_b$, and $d\rho_b(\liez^s)=\bar{\lies}$.  
%%In other words, $\liez^s$ is an abelian $\RR$-split Lie subalgebra of $\End
%%\lieg$, which stabilizes $\iota_b(\lieh)$, and whose restriction  to
%%$\iota_b(\lieh)$, that we denote $(\liez^s)_{|\iota_b(\lieh)}$, is
%%isomorphic to $\bar{\lies}$. 
Then we have
$$\mbox{rk}(\AdgP) \geq \dim \liez^s \geq \dim \bar{\lies} = \rkalg(\Ad H)$$  

\medskip

To prove (2), let $N$ be a connected nilpotent subgroup of $\Ad H$ of maximal degree, and let $\bar{N}$ be the Zariski closure. Let $S = \mbox{Ad}^{-1}(\bar{N}) \subset H$.  As above, $S$ is amenable, so theorem \ref{main-theorem0} gives an algebraic subgroup $\check{S}^b$ such
that $\rho_b(\check{S}^b) = \Sadd$.  The right-hand side is the discompact radical of $\bar{N}$.

\begin{lemma}
\label{nilpdeg-lemma}
Let $N$ be a connected nilpotent group in $\GL(\RR^k)$.  Let $\Sadd$ be the discompact radical of the Zariski closure $\bar{N}$.  Then $n(\bar{\lies}_d) = n(\lien)$.
\end{lemma}

\begin{Pf}
Let $n = n(\lien)$.  Any $n$-tuple $(u_1, \ldots, u_n) \in \lien^n$ satisfies 
$$[u_n, \cdots [u_2,u_1]\cdots] = 0$$
 Because this is an algebraic condition, it is satisfied by all $n$-tuples of elements of $\bar{\lien}$.  Therefore $n(\lien) = n(\bar{\lien})$. 
 
Of course $n(\bar{\lies}_d) \leq n(\bar{\lien})$, so it suffices to show
$n(\bar{\lies}_d) \geq n(\bar{\lien})$.  Write $\bar{\lien}= \lier \ltimes
\lieu$ where $\lieu$ is composed of nilpotents, and $\lier$ is reductive.  The algebra $\lier$ is reductive and nilpotent, so it is abelian. The adjoint
representation of $\lier$ on $\lieu$ is both nilpotent and reductive, so it is
trivial.   Therefore $\lier$ is abelian and centralizes $\lieu$, so that $n(\bar{\lien})=n(\lieu)$. 
But $\lieu \subset \bar{\lies}_d$, so $n(\bar{\lien}) \leq n(\bar{\lies}_d)$.
\end{Pf}  

Now as above $d\rho_b$ is a surjective morphism of algebraic Lie algebras
$\check{\lies}^b \to \bar{\lies}_d$.  Take an algebraic Levi decomposition 
$$\check{\lies}^b = \lier \ltimes \lieu$$ 

 and write $\lier = \liez \oplus \liel$ as above.  Since $\bar{\lies}_d$ is
 nilpotent, and $\liel$ is semisimple, again $\liel \subset \ker d\rho_b$. The kernel of $d\rho_b$ in $\lieu$ is $\ad \lier$-invariant, so there is an
$\ad \lier$-invariant complementary subspace $V$ in $\lieu$.  
The subspace $\liez \oplus V$ is mapped surjectively on $\bar{\lies}_d$ by $d\rho_b$. Let $\liez^{\prime}$ be the kernel in $\liez$ of the adjoint restricted to
$V$, and let $\liez^{\prime \prime}$ be a complementary subalgebra in
$\liez$.  Then $\liez= \liez^{\prime} \oplus \liez^{\prime \prime}$.  Because
elements of $\liez^{\prime \prime}$ act via the adjoint on $V$ by nontrivial
$\CC$-split endomorphisms, none of them is nilpotent. Therefore $\liez^{\prime
  \prime} \subset \ker d\rho_b$.  Let $\lieu^{\prime}$ be the Lie algebra
generated by $V$.  Then the subalgebra $\liez^{\prime} \oplus \lieu^{\prime}$
of $\check{\lies}^b$ is nilpotent and maps onto $\bar{\lies}_d$ via $d \rho_b$.  By lemma \ref{nilpdeg-lemma}, we conclude:
$$ n(\Ad H) = n(\lien)=n(\bar{\lies}_d) \leq n(\liez^{\prime} \oplus
\lieu^{\prime}) \leq n(\mbox{Ad}_{\lieg}P)$$ $\diamondsuit$

\subsection{Semi-discrete groups: proof of \ref{bound-semidiscrete-corollary}}

Recall that $\Gamma$ is a lattice in an algebraic group $L < \GL(V)$ with $L =
L_d$, and the group $H = \Gamma
\ltimes V$ acts by automorphisms of the compact Cartan geometry $M$.  Note
that $L$ is the Zariski closure of $\Ad H$ on $V = \lieh$, by the Borel
density theorem.  

Let $U = \Mon(\lieh,\lieg)$, and let $ \iota : B \rightarrow U$ be the $H
\times P$-equivariant map defined in section \ref{natural-objects-section}.  
Let $Q < L$ be an amenable, algebraic subgroup with $Q_d = Q$.  Much as in
section \ref{main-proof-section}, we seek a $Q$-fixed point in $U/\AdgP$.

The Zariski closure $L$ commutes with $\AdgP$ on $U$ because $H$ does.
Let $B' = L \times_{\Gamma} B$ be the $L$-space of equivalence classes
of pairs $(h,b) \in L \times B$, where $(h,b) \sim (h \gamma^{-1},\gamma.b)$.
Define an $L \times P$-equivariant extension of $\iota$
\begin{eqnarray*}
\iota' & : & B' \rightarrow U\\
 &  & [(h,b)] \mapsto h.\iota(b)
 \end{eqnarray*} 

Let 
$$M' = L \times_{\Gamma} M = B'/P$$
Note that $M'$ is an $M$-bundle over $L/\Gamma$, which has an $\bar{H}$-invariant probability measure
$\nu$.  The space of probability measures on $M'$
projecting to $\nu$ is a compact convex $L$-space, so it contains a
fixed point for $Q$; in other words, there is a finite $Q$-invariant measure
$\mu$ on $M'$.

Denote by $\eta$ the $H$-equivariant map $M \rightarrow U/\AdgP$ covered by
$\iota$. There is an $L$-equivariant continuous lift
$$ \eta' : M' \rightarrow U/\AdgP$$

The $L$-action on $U$ preserves the stratification as in the proof of \ref{main-theorem} in section \ref{main-proof-section}, so
it acts on the smooth quotient varieties $V_i = U_i/\AdgP$, for each stratum
$U_i \subset U$.  Let $\mu_i$ be the restriction of $\eta'_*(\mu)$ to $V_i$.
Each $\mu_i$ is $Q$-invariant.  The corollary \ref{cor:exist-fixed-point} of
the Borel density theorem gives a $Q$-fixed point $v =
\eta'([(h,x)]) = h.\eta(x)$ for some $h \in L$ and $x \in M$.  Then
$hQh^{-1}$ fixes $\eta(x)$.  In other words, for any $b \in B$ lying over $x
\in M$, the algebraic subgroup
$\check{P}^b < \AdgP$, when restricted to $\iota_b(V)$, contains a subgroup isomorphic to $hQh^{-1}$.
Let $\check{S}^b$ be the algebraic subgroup $\rho_b^{-1}(hQh^{-1})$, so
$\rho_b$ is an algebraic epimorphism from $\check{S}^b$ onto $hQh^{-1}$.

Now take $Q < L$ to be an $\RR$-split algebraic subgroup of dimension $\rk L$.  Then following the proof of
\ref{bounds-theorem} (1) gives 
$$\rk L = \dim Q \leq \rk(\mbox{Ad}_{\lieg} P)$$

Similarly, there is a connected, nilpotent, algebraic $Q < L$ of degree
$n(L)$, and following the proof of \ref{bounds-theorem} (2) gives $n(L) \leq
n(\mbox{Ad}_{\lieg} P)$.

\section{A rigidity result for parabolic geometries}

\subsection{Parabolic subalgebras, grading, regular connections}
\label{parabolic-background-subsection}
In this section, $G$ is a connected semisimple Lie group.  Let $\liea  \subset
\lieg$ be a maximal $\RR$-split abelian subalgebra, and $\Pi \subset \liea^*$
the root system. This set of roots admits a basis $\Phi=\{
\alpha_1,...,\alpha_r \}$ of \emph{simple roots}, such that any root can be
written uniquely as a sum
$n_1\alpha_1+...+n_r \alpha_r$, where $n_1,...,n_r$ are integers of the same
sign. The simple roots determine a decomposition 
$$\Pi = \Pi^+ \cup \{  0 \} \cup \Pi^-$$
 where $\Pi^+$ (resp. $\Pi^-$) denotes the positive (resp. negative)
 roots---those for which all the integers $n_i$ are positive (resp. negative).
 
For $\alpha \in \Pi$, the root space
$$\lieg_{\alpha}=\{Y \in \lieg \ | \ (\ad X)(Y)=\alpha(X)\cdot Y \qquad
\forall X \in \liea   \}$$
The associated root space decomposition is
$$\lieg=\Sigma_{\alpha \in \Pi^-}\lieg_{\alpha}\oplus \liem \oplus \liea
\oplus \Sigma_{\alpha \in \Pi^+}\lieg_{\alpha}$$
  The subalgebra $\liem$ is a compact subalgebra, and $\liea \oplus \liem $ is
  exactly the centralizer of $\liea$ in $\lieg$.  We denote by $A$ the
  connected subgroup of $G$ with Lie algebra $\liea$, and by $M_G$ the
  subgroup with Lie algebra $\liem$; it is compact if and only if $Z(G)$ is finite.

Let $\Theta \subset \Phi$, and 
$$\Theta^- = \{ \alpha  = \Sigma n_i \alpha_i \ | \ \alpha_i \in \Phi \
\mbox{and} \ n_i \in \ZZ^{\leq 0} \}$$

The \emph{standard parabolic subalgebra associated to $\Theta$} is
$$\liep = \Sigma_{\alpha \in \Theta^-} \lieg_{\alpha}\oplus \liem \oplus \liea
\oplus \Sigma_{\alpha \in \Pi^+} \lieg_{\alpha}$$ 
For example, taking $\Theta = \emptyset$ yields the {\it minimal parabolic
  subalgebra}
$$\liep=\liem \oplus \liea \oplus \Sigma_{\alpha \in \Pi^+} \lieg_{\alpha}$$
  For $\Theta = \Phi$, one gets $\liep=\lieg$.

Let $\Theta$ be a proper subset of $\Phi$. Reordering the elements of $\Phi$
if necessary, we can assume $\Phi \setminus \Theta = \{
\alpha_1,...,\alpha_m\} $. For every $j \in \ZZ$, define
$$\Pi_j = \{\Sigma_{i=1}^r  n_i\alpha_i \ | \ n_i \in \ZZ^{\leq 0} \ \forall \ i
\ \mbox{or} \ n_i \in \ZZ^{\geq 0} \ \forall\ i \ \mbox{and} \ \Sigma_{i=1}^m
n_i=j \}$$
Then put $\lieg_{j}=\Sigma_{\alpha \in \Pi_j} \lieg_{\alpha}$, and $\lieg^j=\Sigma_{i \geq j} \lieg_i$.

Let $k$ be the greatest integer such that $\Pi_k \not = \emptyset$.
Associated  to the parabolic subalgebra $\liep$ are a grading  
$$\lieg=\lieg_{-k} \oplus \cdots \oplus \lieg_{-1} \oplus \lieg_0 \oplus
\lieg_1 \oplus \cdots \oplus \lieg_k$$  
as well as a filtration 
$$\lieg^{k} \subset \lieg^{k-1} \subset ... \subset \lieg^{-k}=\lieg$$ 
with $[\lieg^i,\lieg^j] \subset \lieg^{i+j}$.  Notice that $\liep=\lieg^0$.

\subsection{Parabolic geometries}
\label{parabolic-geometries}
Let $G$ be as above, and $P < G$ a parabolic subgroup with Lie algebra
$\liep$.  Let $\lieg_{i}$, for $i = -k, \ldots, k$, be the subspaces of the
associated grading, as above.  A {\it parabolic geometry} is a Cartan geometry $(M,B,\omega)$ modeled on $\XX=G/P$. As illustrated in section \ref{intro-parabolic-section} of the introduction, many interesting examples of geometric structures have an interpretation as a parabolic geometry.

The Cartan connection $\omega$ of a parabolic geometry modeled on $\XX=G/P$ is
{\it regular} if the corresponding curvature function satisfies for any $i,j <0$
$$\kappa(\lieg^i, \lieg^j) \subset \lieg^{i+j+1}$$
 The regularity assumption on a Cartan connection is natural.  Indeed, in
 every case in which the equivalence problem is solved for parabolic geometries, regularity is one of the normalization conditions on the Cartan connection to make it unique.

%\subsection{Rates of contraction of sequences in a Cartan subgroup}
 %Any other Cartan subalgebra of $\liep$ (i.e any maximal $\RR$-split abelian subalgebra) is conjugated to $\liea$ by an element  $p \in P$. So, it makes sense to define $\liea^+$ as the conjuThe closed group $A=Exp(\liea)$ is called a Cartan subgroup of $P$.  We put any norm $||.||$ on $\liea$. 

\subsection{Flatness and completeness:  proof of \ref{parabolic-flatness-theorem}}
\label{proof-parabolic-flatness-theorem}
Let $G$ be a semisimple Lie group and $P$ a parabolic subgroup as above.
Throughout this section, $(M,B,\omega)$ is a Cartan geometry modeled on
$\XX=G/P$ with $\omega$ regular. 

Inside $\liea$, define the positive Weyl chamber 
$$\liea^+ = \{X \in \liea \ | \ \alpha(X) >0 \ \forall \ \alpha \in \Pi^+ \}$$ 

 Let $(a_k)$ be a sequence of $A$ tending to infinity---that is, leaving every
 compact subset of $A$.  Fix any norm on $\lieg$, and write 
$$a_k=\exp(X_k)=\exp(||X_k||u_k)$$
 where $|| u_k|| = 1$.  An {\it an asymptotic direction of $(a_k)$} is any
 cluster point $u \in \liea$ of the sequence $(u_k)$.

The notion of {\it holonomy sequences} will be important in the proof
below. Let $H < \Aut M$, for $(M,B,\omega)$ a Cartan geometry modeled on
$G/P$, let $b \in B$, and $x = \pi(b)$. For $h_k \in H$ such
that $h_k.x \rightarrow y_{\infty} \in M$, a {\it holonomy sequence associated
  to $(b,h_k)$}, is a sequence $(p_k)$ of $P$ such that 
$$h_k.b.p_k^{-1} \rightarrow b_{\infty} \in \pi^{-1}(y_{\infty})$$
 If $(p_k)$ is a holonomy sequence associated to $(b,h_k)$, then all other
 holonomy sequences associated to $(b,h_k)$ are $(l_kp_k)$ with $(l_k)$ a
 convergent sequence of $P$.  Also, if $b^{\prime}=b.p^{\prime}$, then any
 holonomy sequence associated to $(b^{\prime},h_k)$ is of the form
 $(p_kp^{\prime})$, where $(p_k)$ is a holonomy sequence associated to
 $(b,h_k)$. These observations lead naturally to the following notion of equivalence.

\begin{definition}
%\label{equivalent}
Two sequences $(p_k)$ and $(p_k^{\prime})$ are \emph{equivalent} if there are a convergent sequence $(l_k)$ in $P$ and $q \in P$ such that $p_k^{\prime}=l_kp_kq$, for all $k \in \NN$.
\end{definition}

Note that a sequence of $P$ and its conjugate by some $p \in P$ are equivalent. Also, two converging sequences of $P$ are equivalent.

\begin{definition}
\label{equivalent}
A \emph{holonomy sequence at $b \in B$} is a holonomy sequence associated to
$(b,h_k)$ for some $(h_k)$ in $H$.  For $x \in M$, a \emph{holonomy sequence
  at $x$} is a holonomy sequence at some $b \in B$ for which $\pi(b) = x$. 
\end{definition}

Remark that if $(p_k)$ is a holonomy sequence at $x$, then any sequence equivalent to $(p_k)$ is also a holonomy sequence at $x$.

If $b \in B$ is such that $\kappa_b=0$, then $\kappa_{b^{\prime}}=0$ for every $b^{\prime}$ in the same fiber. In the following, by a slight abuse of language, we will say that the curvature function vanishes at $x \in M$, if $\kappa_b=0$ for all $b$ in the fiber of $x$.
We can now prove

\begin{proposition}
\label{vanishing}
Let $H < \Aut M$ be a Lie subgroup, and $S< H$.
\begin{enumerate}
\item{If there is $b \in B$ such that $\Ad A \subset \check S^b$, then $\kappa_b=0$.} 
\item{Suppose there is $x \in M$ such that for any $u \in \liea^+$, there is a
    holonomy sequence  $(p_k) \subset A$ at $x$ with asymptotic direction $u$.  Then there is an open neighborhood $U$ of $x$ on which the curvature vanishes.}  \end{enumerate}
\end{proposition}

\begin{Pf}
We begin with the proof of point $(2)$.
Let $i \leq j<0$, and $l>i+j$. Let $\alpha \in \Pi_i$, $\beta \in \Pi_{j}$,
and $\nu \in \Pi_{l}$  be three roots of $\lieg$.  Choose nonzero vectors $v
\in \lieg_{\alpha}$ and $w \in \lieg_{\beta}$.  Choose $u \in \liea^+$, such
that $\alpha(u)+\beta(u) < \nu(u)$ (Such $u$ always exits, since $\alpha +
\beta$ is in $\Pi_{i+j}$, hence cannot equal $\nu$).  Finally, choose $(p_k)$
a holonomy sequence at $x$ with $u$ for asymptotic direction.  Passing to a
subsequence if necessary, write $p_k=\exp(X_k)$ with 
$$\lim_{k \to \infty} \frac{X_k}{||X_k||}=u$$
 By definition of $(p_k)$, there is a sequence $(h_k)$ of $H$ and $b \in B$
 over $x$ such that $h_k.x \rightarrow x_{\infty} \in M$, and
 $$b_k=h_k.b.p_k^{-1} \rightarrow b_{\infty} \in \pi^{-1}(x_{\infty})$$
For all $k \in \NN$
\[ (\Ad p_k).\kappa_{b}(v,w)=\kappa_{b_k}((\Ad p_k).v,(\Ad p_k).w)=e^{\alpha(X_k)+\beta(X_k)} \kappa_{b_k}(v,w)\]

Since $\alpha(u) + \beta(u)-\nu(u)<0$, 
$$e^{-\nu(X_k)}(\Ad p_k).\kappa_{b}(v,w) \rightarrow 0$$   
Then the component of $\kappa_b(v,w)$ on $\lieg_{\nu}$ is zero. Since this is true for any $\nu \in \Pi_{l}$ with $l>i+j$, the component of $\kappa_b(v,w)$ on $\lieg^{i+j+1}$ is zero, which implies $\kappa_b(v,w)=0$ by the regularity condition. Then at any point $x \in M$ where the assumptions of point $(2)$ hold, the curvature vanishes.

Let $b \in B$ and $\xi \in T_bB$. Let $X_{\xi}$ be the unique vector field on $B$ having the  property that $\omega_{b^{\prime}}(X_{\xi}(b^{\prime}))=\omega_b(\xi)$ for all $b^{\prime} \in B$. If $\phi^t$ is the flow generated by $X_{\xi}$, and if $\phi^1$ is defined, then put $\exp_b(\xi)=\phi^1.b$. The map $\exp_b$ is defined on a neighborhood of $0_b$ in $T_bB$, and yields a diffeomorphism of a sufficiently small neighborhood of $0_b$ onto its image. 

Let $u$ be in $\liea^+$ and $(p_k)$ a holonomy sequence at $x$ with asymptotic
direction $u$.  Let $h_k \in H$ and $b \in \pi^{-1}(x)$ such that 
$$b_k=h_k.b.p_k^{-1} \rightarrow b_{\infty} \in B$$
 Choose an open $U \subset \lien^-=\lieg_{-k}\oplus \cdots \oplus \lieg_{-1}$ containing $0_{\lien^-}$ such that $\exp_b$ is a diffeomorphism from $\omega_b^{-1}(U)$ onto its image; note that $\pi \circ \exp_b(\omega_b^{-1}(U))$ is open in $M$.  For every $k \in \NN$,

\[h_k.\exp_b(\omega_b^{-1}(U)).p_k^{-1} = \exp_{b_k}(\omega_{b_k}^{-1}((\Ad p_k).U)) \]

Taking a subsequence if necessary, we assume that $u$ is the only asymptotic
direction of  $(p_k)$. Since $u \in \liea^+$, the elements $\Ad p_k$ act by
contraction on $\lien^-$, which implies that $(\Ad p_k).U \rightarrow 0_{\lieg}$
for the Hausdorff topology. As a consequence, 
$$\exp_{b_k}(\omega_{b_k}^{-1}((\Ad p_k).U)) =
h_k.\exp_b(\omega_b^{-1}(U)).p_k^{-1} \rightarrow \{b_{\infty}\}$$ 
Let $V=\pi \circ \exp_b(\omega_b^{-1}(U))$.  For any $x^{\prime} \in V$, the
sequence $(p_k)$ is a holonomy sequence at $x^{\prime}$. The hypotheses of point $(2)$ of the proposition are thus satisfied at any $x^{\prime} \in V$. By what was said above, the curvature vanishes at every point of the open set $V$.

We now prove point $(1)$. Since $\check S^b$ preserves $\kappa_b$, any
sequence $(p_k)$ of $A$ satisfies
$$(\Ad p_k).\kappa_{b}(v,w)=\kappa_{b}((\Ad p_k).v,(\Ad p_k).w) \ \forall \ v,w \in
\lieg$$
 The proof is then basically the same as the first part of that of point $(2)$.
\end{Pf}

We can now prove the first point of theorem
\ref{parabolic-flatness-theorem}. By the hypothesis on $\rkalg(\Ad H)$, there
is an abelian $\RR$-split subgroup of $ \Ad H$ with Zariski closure $\bar S$
of dimension $\rk G$.  Let $S= \Ad^{-1}(\bar S) \subset H$. Since the kernel
of $\Ad$ on $H$ is amenable by assumption (\S \ref{intro-mainthm-section}),   
$S$ is amenable as well.  Theorem \ref{main-theorem} gives $b \in B$ such
that $\rho_b$ maps $\check S^b$ surjectively on $\bar S=\Sadd$.  

As in the proof of theorem \ref{bounds-theorem} in section \ref{bounds-proof-section}, there is an abelian
$\RR$-split subalgebra  $\liez^s \subset \check \lies^b$ such that $d
\rho_b(\liez^s)=\bar \lies$. Then $\dim \liez^s = \rk G$, and there exists a
maximal $\RR$-split subalgebra $\liea^{\prime} \subset \liep$ such that $\ad
\liea^{\prime}=\liez^s$.  There exists $p_o \in P$ such that $(\Ad
p_0)(\liea^{\prime}) = \liea$. Since 
$$\Ad A < \check S^{b.p_0}=(\Ad p_0)\check S^b(\Ad p_0^{-1})$$
 we can apply point $(1)$ of proposition \ref{vanishing}. Together with lemma \ref{curvature-lemma}, it yields that $-\iota_{b.p_0} : \lieh \to \lieg$ is an embedding of Lie algebras, proving the first point of the theorem.

We now come to the second  part of theorem
\ref{parabolic-flatness-theorem}. By assumption, there is a connected Lie group $A_H < H$ with Lie algebra
$\liea_H \subset \lieh$, such that $\Ad A_H$ is  $\RR$-split, abelian, and has
dimension $\rk G$. Moreover, since the kernel of $\mbox{Ad}$ is amenable by
assumption, so is $A_H$. %connected subgroup now assumed to be semisimple. Let $\liea_H$ be a Cartan subalgebra of $\lieh$, and $\Pi_H$ an associated set of roots. The Lie algebra $\lieh$ writes as $\lieh = \Sigma_{\beta \in \Pi_H^+}\lieh_{-\beta}\oplus \liem_H \oplus \liea_H\oplus \Sigma_{\beta \in \Pi_H^+}\lieh_{\beta}$, where $\liem_H$ is a compact Lie algebra, and $\liea_H \oplus \liem_H$ is the centralizer of $\liea_H$. 
We first show
\begin{lemma}
\label{fix}
Every closed $A_H$-invariant subset of $M$ contains an $A_H$-fixed point
$x$ with an open neighborhood on which the curvature vanishes.
\end{lemma}

\begin{Pf}
Let $F$ be a closed $A_H$-invariant subset.  Apply theorem \ref{main-theorem}
with $S=A_H$.  It gives $b \in \pi^{-1}(F)$ such that $\check S^b$ is mapped
surjectively by $\rho_b$ on $(\bar{A}_H)_d$.  Exactly as in the proof of the first point, changing $b$ into $b.p_0$ if necessary, we obtain $\Ad A < \check S^b$, and $\rho_b$ maps $\Ad A$ onto the identity component of $(\bar{A}_H)_d$.  
%Since  $ \Ad A_H$ has already the maximal dimension allowed for an $\RR$-split connected abelian subgroup of $G$, we get $(\bar{A}_H)_d= \Ad A_H$, and $\rho_b$ is an isomorphism between $A$ and $\Ad A_H$. 

Now $(\bar{A}_H)_d$ acts trivially on $\liea_H$, so $\Ad A$ acts trivially on
$\liea_H^b=\iota_b(\liea_H)$. Thus $\liea_H^b \subset \liea \oplus \liem$, the
centralizer of $\liea$; in particular, $\liea_H^b \subset \liep$, which proves
that $A_H$ fixes $x=\pi(b)$.  

Since $\liea_H^b \subset \liep$, lemma \ref{curvature-lemma} ensures that
$-\iota_b$ is a Lie algebra isomorphism between $\liea_H$ and $\liea_H^b$, so
$\liea_H^b$ is abelian. In fact, for $u \in \liea_H$ and $v \in \lieh$
$$(\ad \iota_b(u))(\iota_b(v))=-\iota_b((\ad u).v)$$
 In particular, $\ad \liea_H^b$ acts on $\iota_b(\lieh)$ as an $\RR$-split
 subalgebra of dimension $\rk G$.  Decompose $u \in \liea_H^b$ into
 $u_{\liea}+u_{\liem}$.  Because $\liea$ and $\liem$ are simultaneously
 diagonalizable on $\lieg$ and $u$ preserves $\iota_b(\lieh)$, both $\ad u_{\liea}$ and $\ad u_{\liem}$ preserve
 $\iota_b(\lieh)$; moreover, the restriction of $\ad u_{\liem}$ is
 trivial. Then the projection of $\liea_H^b$ modulo $\liem$ is a surjection
 onto $\liea$. 

Let $A_H^b$ be the connected subgroup of $P$ with Lie algebra $\liea_H^b$. Any
sequence in $A_H^b$ is a holonomy sequence at the fixed point $x=\pi(b)$.
From above, $A_H^b < A \times M$, and projects onto the first
factor. Thus any sequence of $A$ is equivalent to  a sequence in $A_H^b$, and
it follows from the remark after definition \ref{equivalent} that any sequence
of $A$ is a holonomy sequence at $x$. We can now apply proposition
\ref{vanishing} to obtain an open neighborhood $V$ containing $x$ on which the curvature vanishes. \end{Pf}

\begin{corollary}
The Cartan geometry $(M,B,\omega)$ is flat.
\end{corollary}

\begin{Pf}
Let $N_0 \subset M$ be the subset on which the curvature of $\omega$
vanishes. By the previous lemma, $N_0$ has nonempty interior. If $x \in
\partial N_0$, the topological boundary of $N_0$ in $M$, then
$\overline{H.x}$, the closure of the $H$-orbit of $x$, is closed and
$H$-invariant. By lemma \ref{fix}, there are $y \in \overline{H.x}$ and an
open neighborhood $V$ of $y$ such that the curvature vanishes on $V$. Then $y$
could not have been in $\partial N_0$, so $\partial N_0=\emptyset$, which means $N_0=M$.   \end{Pf}

\begin{corollary}
\label{complete-corollary}
The geometry $(M,B,\omega)$ is geometrically isomorphic to $\Gamma \backslash \widetilde \XX$. 
\end{corollary}

\begin{Pf}
By the corollary directly above, the geometry $(M,B,\omega)$ is flat---in
other words, there is a $(G,\XX)$-structure on $M$. The proof of the
completeness of this structure is a straightforward generalization of that in
\cite{fz}. The $(G,\XX)$-structure on $M$ defines a developing map $\delta :
\widetilde M \to \XX$, which is a geometric immersion, and a holonomy morphism
$\rho : \Aut \widetilde{M} \to G$ satisfying the equivariance relation $\rho
\circ \delta=\delta \circ \rho$ (see \cite{thurston}).  Let $\pi_1(M)$ be the
fundamental group of $M$, which acts by automorphisms on $\widetilde M$, and consider $\Hol(M)=\rho(\pi_1(M)) \subset G$,  the {\it holonomy group} of $M$.

At the infinitesimal level, there is a morphism of Lie algebras 
$$d \rho_e : \chi^{kill}(\widetilde M) \to \lieg$$
 where $\chi^{kill}(\widetilde M)$ denotes the Lie algebra of Killing fields
 on $\widetilde M$.  Since $H$ acts by automorphisms on $M$, there is an
 infinitesimal action of $\lieh$ by Killing fields on $\widetilde M$, so that
 $\lieh \hookrightarrow \chi^{kill}(\widetilde M)$.  Observe that any Killing
 field of $\widetilde M$ lifted from $M$ is centralized by $\pi_1(M)$. 

Now the morphism $d \rho_e : \lieh \to \lieg$ maps $\lieh$
(resp. $\liea_H$) injectively onto a  Lie subalgebra $\lieh^{\prime} \subset
\lieg$ (resp. $\liea_H^{\prime} \subset \lieh^{\prime}$), which moreover is
centralized by $\Hol(M)$. It was established in the proof of lemma \ref{fix}
that $\liea_H^{\prime}$ is included, up to conjugacy, in the sum $\liea \oplus
\liem$, and projects onto $\liea$.  Choose $u_1,\ldots,u_s$ spanning $\liea_H^{\prime}$, with projections spanning $\liea$. Since those projections
are just the $\RR$-split parts of $u_1,\ldots,u_s$,
they are centralized by $\Hol(M)$. Then $\Hol(M)$ centralizes $\liea$, so is
contained in $A \times M_G$.  Now $\Hol(M)$ centralizes $\lieh^{\prime}$, and
$A$ acts faithfully on $\lieh^{\prime}$ because its action is conjugated to
that of $\Ad A_H$ on $\lieh$. It follows that the  projection of $\Hol(M)$ on
$A$ is trivial, namely $\Hol(M) \subset M_G$. The compact group $M_G$ is included in $K$, a maximal compact subgroup of $G$, which acts transitively on $\XX$. 

In conclusion, the $(G,\XX)$-structure on $M$ can be seen as a
$(K,\XX)$-structure. These latter structures are known to be \emph{complete} on
closed manifolds (see proposition 3.4.10 of \cite{thurston}): the developing
map $\delta$ is a covering map (The key point is that $\delta$ is not changed
by viewing the Cartan geometry as a $(K,\XX)$-structure rather than a
$(G,\XX)$-structure). Completeness means $\widetilde{M}$ is geometrically
isomorphic to  $\widetilde \XX$, so $\Aut \widetilde{M} \cong \widetilde{G}$,
and $M$ is geometrically isomorphic to $\Gamma \backslash \widetilde{\XX}$,
with $\Gamma$ a discrete subgroup of $\widetilde{G}$ acting properly
discontinuously on $\widetilde{\XX}$.
\end{Pf}

\section{Results about local freeness}

%Results linked with the question of local freeness for volume-preserving actions were already proved in  \cite{zeghib1}, \cite{stuck-zimmer}, \cite{}.

\subsection{Local freeness almost everywhere: proof of \ref{local-freeness-theorem}}

Here $(M,B,\omega)$ is a connected Cartan geometry modeled on $G/P$, and $H
< \Aut M$ preserves a finite measure $\nu$ on $M$ with full support.  Theorem \ref{main-theorem} yields a subset $\Omega \subset M$ of
full measure such that for every $b \in B$ projecting on $x \in \Omega$, the
image $\rho_b(\check{H}^b)=\Hadd$.  The stabilizer Lie algebra $\lieh_x$ has
image $\iota_b(\lieh_x) = \iota_b(\lieh) \cap \liep$.  Since
$\iota_b(\lieh) \cap \liep$ is invariant by $ \check{H}^b$, the Lie algebra $\lieh_x$ is invariant by $\Hadd$, as desired.

%hence by $\Ad H$. It comes that the identity component is normal in $H$. 

Assume now that $\Hadd$ leaves stable only a finite family $\{ 0 \} = F_0$,\ldots,$F_s$ of
subalgebras in $\lieh$. Write $\Omega=\Omega_0
\cup \cdots \cup \Omega_s$, where $\Omega_i$ contains all $x \in \Omega$ such that
$\lieh_x=F_i$.   Let $H_i$ be the connected Lie subgroup of $H$ with Lie
algebra $F_i$. The action of $H_i$ is faithful, so that for $i \not = 0$, the
interior of ${\Omega}_i$ is empty. This follows from the lemma

\begin{lemma}
Let $\{ h^t \}_{t \in \RR} < \Aut M$ be a $1$-parameter subgroup, and $U \subset M$ an open subset such that $\left. h^t
\right|_U=Id$ for every $t \in \RR$. Then $h^t = Id$ for every $t \in \RR$.

\end{lemma} 
\begin{Pf}
As already noticed, the action of $h^t$ on $B$ is free because it preserves a
framing. So if suffices to show that $h^t$ fixes a point on $B$ for every $t
\in \RR$. Pick $x \in U$ and $b \in B$ projecting on $x$.  Thanks to the
Cartan connection, it is  possible to define a developing map ${\cal D}_x^b$
from the space $C^1([- 1,1],U)$ of $C^1$ curves $\gamma$ in $U$ with
$\gamma(0)=x$ to the space $C^1([-1,1],\XX)$ of $C^1$ curves $\beta$ in $\XX$
with $\beta(0)=o$, the class of $P$ in $\XX=G/P$. We refer to \cite{Sharpe}
for the construction of the developing map. Now there is a $1$-parameter
subgroup $\{ p^t \} <  P$ such that $h^t.b.p^{-t}=b$. It is not difficult to
show that ${\cal D}_x^b(h^t.\gamma)= p^t.{\cal D}_x^b(\gamma)$. Because
$\left. h^t\right|_U = Id$, all developments of curves in $U$ are pointwise
fixed by $\Ad p^t$.  The union of all these developments is a subset
of $\XX$ with nonempty interior, pointwise fixed by $p^t$. As a consequence, $p^t$ acts trivially on $\XX$ for every $t \in \RR$. Since the kernel of the action of $P$ on $\XX$ was assumed to be finite,  $p^t=e$ for every $t \in \RR$.
\end{Pf}

Now because $\Omega$ is conull for a measure of full support, it is dense in $M$, and $\Omega_1 \cup \cdots \cup \Omega_s$ has
empty interior. Thus $\Omega_0$, on which $H$ is locally free, is dense in
$M$. Since the subset of $M$ where the action is locally free is clearly open,
theorem \ref{local-freeness-theorem} follows. $\diamondsuit$ 

\subsection{Open orbit implies homogeneity: proof of \ref{homogeneous-theorem}}
In this section, we assume that $H < \Aut M$ is a Lie subgroup such that $\Had$, the Zariski closure
 of $\Ad H$ in $\Aut \lieh$, has no algebraic compact quotients---that is, $\Hadd = \Had$. The Cartan geometry $(M,B,\omega)$ is connected and $n$-dimensional, equipped with an $H$-invariant volume form $\nu$ of finite total volume.
 
Recall that $\sigma$ is the ($\AdgP$)-equivariant projection from $\lieg$ onto $\lieg/\liep$. 
 For every $b \in B$ projecting on $x \in M$, there is a natural isomorphism
 $j_b : T_xM \to \lieg/\liep$ with $j_b(u)=\sigma \circ \omega_b(\check u)$,
 where $\check u \in T_b B$ projects on $u$. It is easy to check that this
 isomorphism is well-defined, and $j_{b.p}=(\Ad p^{-1})\circ j_b$ for every $p \in P$ and $b \in B$.

By this remark, we see that a vector field on $M$ is nothing else than a map
$\phi : B \to \lieg/\liep$ satisfying $\phi(b.p)=(\Ad p^{-1}). \phi(b)$. In the
same way, a continuous $n$-form $\nu$ on the manifold $M$ can be seen as a
continuous map  $\nu : B \to \wedge^n(\lieg/\liep)^*$, such that for any
$n$-tuple $(u_1,\ldots,u_n)$ of $(\lieg/\liep)^n$, 
$$\nu(b.p)(u_1,\ldots,u_n)=\nu(b)((\Ad p).u_1,\ldots,(\Ad p).u_n)$$
 Our hypothesis is that $\nu(h.b)=\nu(b)$ for all $h \in H$. 

Let $U$ be the real algebraic variety $\Mon(\lieh,\lieg) \times
\wedge^n(\lieg/\liep)^*$. The map $\phi : b \mapsto (\iota_b, \nu(b))$ is
continuous and $H \times P$-equivariant. Here, the action of $H$ on
$\wedge^n(\lieg/\liep)^*$ is trivial, and for $\alpha \in
\wedge^n(\lieg/\liep)^*$ and $\check p \in \AdgP$,  
$$(\check p.\alpha)(u_1,\ldots,u_n)=\alpha(\check p.u_1,\ldots,\check p.u_n).$$

For every $b \in B$, we define, as in section \ref{natural-objects-section}

\[ \check{H}^b = \{\check p \in \check{P}^b \ | \ \rho_b(\check p) \in \Hadd, \ \mbox{and} \  \check p.\nu(b)=\nu(b) \} \]

Following the same proof as for theorem \ref{main-theorem}, we obtain a
$P$-invariant subset $\Lambda \subset B$ projecting on a set of full
$\nu$-measure, such that for every $b \in \Lambda$, $\rho_b$ is a surjection
from $\check H^b$ onto $\Had=\Hadd$. In particular, if $b \in \Lambda$ and $h
\in H$, then there is $\check p_h \in \check H^b$, not necessarily unique, such that $\iota_{h.b}= \check p_h. \iota_b$ and $\check p_h.\nu(b)=\nu(b)$.

%As before, we denote by $\Mon(\lieh,\lieg)$ the space of injective linear maps
%from $\lieh$ to $\lieg$. 
%% The Cartan geometry on $M$ is supposed to be reductive, so that
%% $\lieg=\lien \oplus \liep$, where $\lien$ is $\AdgP$-invariant. There is a
%% natural projection $\pi : \lieg \to \lien$, with kernel $\liep$. This
%% projection is $\AdgP$-equivariant. We introduce $\Mon(\lieh,\lieg)^T= \{
%% \alpha \in \Mon(\lieh,\lieg) \ ���� \ \pi_{����\alpha(\lieh)} \ is \
%% surjective\}$.    

 Observe that if $M$ has an open $H$-orbit, then it has an open orbit for
 $H^o$, the identity component of $H$. Indeed, if $H^o.x$ is not open for some
 $x \in M$, then it has volume zero. But $H$ is a Lie group, so it has
 countably-many connected components. Thus $H.x$ is the union of
 countably-many sets of volume zero, so that it also has volume zero,
 contradicting the assumption that it was open. 

Let $\cal{O}$ be an open orbit of $H^o$ on $M$. Since any set of full measure
for $\nu$ intersects $\cal{O}$, there is $b \in \Lambda$ projecting on $x \in \cal{O}$.  Denote by $\overline {\cal{O}}$ the topological closure of $\cal{O}$ in $M$. Pick $x_{\infty} \in \overline {\cal{O}}$ and $b_{\infty} \in B$ projecting on $x_{\infty}$. 
 There is $(h_k)$ a sequence of $H^o$ such
that $h_k.x \rightarrow x_{\infty}$, and thus a sequence $(p_k)$ of $P$ such
that $h_k.b.p_k^{-1} \rightarrow b_{\infty}$.  Now, the tangent space
$T_x(H^o.x) = T_xM$, so there are $X_1,\ldots,X_n \in \lieh$ such that $(\sigma \circ \iota_b(X_1), \ldots, \sigma \circ \iota_b(X_n))$ is a basis of $\lieg/\liep$, and thus $\nu(b)(\sigma \circ \iota_b(X_1), \ldots, \sigma \circ \iota_b(X_n)) \not = 0$.

The equivariance property implies that for every $p \in P$, 
$$\nu(b.p)(\iota_{b.p}(X_1),\ldots,\iota_{b.p}(X_n))=\nu(b)(X_1,\ldots,X_n)$$
 The above equality together with $H$-invariance of $\nu$ implies that there are $\check g_k \in \check H^b$ such that
\begin{eqnarray*}
\nu(h_k.b.p_k^{-1})(\sigma \circ \iota_{h_k.b.p_k^{-1}}(X_1),\ldots,\sigma \circ \iota_{h_k.b.p_k^{-1}}(X_n))
& \\
=\nu(b)(\check g_k.\sigma(X_1),\ldots, \check g_k.\sigma(X_n)) &  
\end{eqnarray*}
 Since $\check g_k.\nu(b)=\nu(b)$, the right-hand side equals $\nu(b)(\sigma
 \circ \iota_b(X_1), \ldots, \sigma \circ \iota_b(X_n))$, and the left one
 tends to $\nu(b_{\infty})(\sigma \circ \iota_{b_{\infty}}(X_1), \ldots,
 \sigma \circ \iota_{b_{\infty}}(X_n))$.  This last term is thus nonzero,
 which proves that $\sigma \circ \iota_{b_{\infty}}(\lieh)=\lieg/\liep$. It
 follows that the $H^o$-orbit of $x_{\infty}$ is open, and that $x_{\infty}
 \in \cal{O}$. Because $M$ is connected, $M=\cal{O}$, which proves that $H^o$ acts transitively on $M$.

Theorem \ref{local-freeness-theorem} implies that for almost every $x \in M$,
the stabilizer subalgebra $\lieh_x$ is $\Had$-invariant. Because $H^0$ acts
transitively, the stabilizer subalgebra is the same at every point of
$M$. Since $H$ acts faithfully, $\lieh_x=\{ 0 \}$ for every $x \in M$.  Then
$M=H^o/\Gamma$, where $\Gamma$ is a discrete subgroup of $H^o$. The $n$-form
$\nu$ is left-invariant, so it induces Haar measure on $H^o/\Gamma$.  Since $M$ has finite volume, $\Gamma$ must be a lattice in $H^o$.  

%% Finally, given that the action of $(\Aut M)^o$ is transitive on $M$, it has
%% dimension at most $n$. But in fact, $\text{dim } (\Aut M)^o=n=\text{dim } H^o$, because the above arguments applied to $(\Aut M)^o$ prove that its action is locally free. This yields $H^o= (\Aut M)^o$.
 
{\bf Deduction of Corollary \ref{gromov-corollary}}

Let $\Autloc (M,\omega,\nu)$ (resp. $\Autloc (M, \cal{S},\nu)$) be the
pseudo-group of local automorphisms of the Cartan geometry (resp. of the
structure $\cal{S}$) that also preserve $\nu$. The hypotheses of corollary
\ref{gromov-corollary} imply $\Autloc (M,\omega,\nu) = \Autloc (M,
\cal{S},\nu)$. The fact that $H$ has a dense orbit ensures that $\Autloc
(M,\omega, \nu)$, hence also $\Autloc (M, \cal{S}, \nu)$, has a dense orbit. 

Adding an analytic volume form $\nu$ to the analytic rigid geometric structure
of algebraic type $\cal{S}$ still yields an analytic rigid geometric structure
of algebraic type $\cal{S}^{\prime}=(\cal{S},\nu)$. Gromov's open-dense
theorem (\cite{Gromov} 3.3.A), implies that $\Autloc (M,\cal{S}, \nu)$ has an open
dense orbit $U$. Note that the local Killing fields for this structure also act
transitively on $U$.

Because the structure on $M$ is analytic and $M$ is simply connected, local
Killing fields extend to all of $M$ (\cite{Amores}, \cite{Nomizu}, \cite{Gromov}). Because $M$ is compact the extended
Killing fields are complete.  Then $\Aut (M,\cal{S},\nu)$ has an open orbit on
$M$, hence so does $\Aut (M,\omega,\nu)$.  The corollary now follows from theorem \ref{homogeneous-theorem}.  $\diamondsuit$

 \subsection{Homogeneous reductive Cartan geometries}
%% It remains to understand the homogeneous case, proving theorem
%% \ref{homogeneous-2-theorem}. 

Finally, we prove theorem \ref{homogeneous-2-theorem}.  The assumptions here are that the geometry is reductive, and the $H$-action is
transitive and faithful, and preserves a finite volume.  By theorem
\ref{homogeneous-theorem}, the action is everywhere locally free.  

Now use
Theorem \ref{main-theorem} to obtain $b \in B$ such that $\rho_b$ is a
surjection from $\check H^b$ to $\Had$. Let $\check H = \mbox{Ad}_{\lieg}^{-1}(\check H^b) <
P$. Because $H$ acts locally freely,  $\iota_b(\lieh)$ is transverse to
$\liep$ for all $b \in B$. Let $\sigma$ be the $\AdgP$-equivariant projection
from $\lieg$ to $\lien$ associated to the decomposition $\lieg=\lien \oplus
\liep$. If $\Ad p$ is trivial on $\iota_b(\lieh)$, then it is trivial on $\lien$. Because $\AdgP$ is faithful on $\lien$, the kernel of
$\rho = \rho_b \circ \mbox{Ad}_{\lieg} :\check H \to \Had$ is trivial, so $\rho$ is an
isomorphism. The map $\iota=\sigma \circ \iota_b$ is an isomorphism from
$\lieh$ to $\lien$, and by theorem \ref{main-theorem}, it intertwines the
representations of $\Had$ on $\lieh$ and $\Ad \check H$ on $\lien$. Set $H^{\prime}=\rho^{-1}(\Ad H)$. When $H$ is connected, then $H$ is a central extension of $H^{\prime}$ as claimed in the second point of theorem \ref{homogeneous-2-theorem}. $\diamondsuit$

\end{document}